\documentclass{eage2022}

\usepackage{amsmath,amssymb,bm}
\usepackage{natbib}
\usepackage{graphicx}
\usepackage[colorlinks]{hyperref}
\usepackage{hologo}
\usepackage{caption}
\usepackage{subcaption}

\usepackage{algpseudocode,algorithm,algorithmicx}
\newcommand*\Let[2]{\State #1 $\gets$ #2}
\algrenewcommand\algorithmicrequire{\textbf{Input:}}
\algrenewcommand\algorithmicensure{\textbf{Output:}}

\usepackage{mathtools}
\renewcommand{\O}{\mathcal{O}}
\newcommand{\bmm}{\mathbf{m}}
\newcommand{\bx}{\mathbf{x}}
\newcommand{\bxs}{\bx_{\mathrm{s}}}

\newcommand{\bu}{\mathbf{u}}
\newcommand{\bd}{\mathbf{d}}
\newcommand{\bnd}{\mathbf{n}_{\mathrm{d}}}
\newcommand{\bnq}{\mathbf{n}_{\mathrm{q}}}
\newcommand{\bq}{\mathbf{q}}
\newcommand{\by}{\mathbf{y}}
\newcommand{\br}{\mathbf{r}}
\newcommand{\brm}{\br(\bmm)}
\newcommand{\F}{\mathrm{F}}
\newcommand{\Fm}{\F(\bmm)}
\newcommand{\A}{\mathrm{A}}
\newcommand{\Am}{\A(\bmm)}
\newcommand{\R}{\mathrm{R}}
\newcommand{\Jac}{\mathrm{J}}
\newcommand{\calN}{\mathcal{N}}
\newcommand{\calL}{\mathcal{L}}
\newcommand{\norm}[1]{\left\lVert#1\right\rVert}
\newcommand{\Sd}{\Sigma_{\mathrm{d}}}
\newcommand{\tSd}{\widetilde{\Sigma}_{\mathrm{d}}}
\newcommand{\Sq}{\Sigma_{\mathrm{q}}}
\newcommand{\<}{\langle}
\renewcommand{\>}{\rangle}
\newcommand{\defeq}{:=}
\newcommand{\eqdef}{=:}
\newcommand{\inv}{^{-1}}
\newcommand{\adj}{^*}
\newcommand{\mean}{\mathbb{E}}
\newcommand{\bz}{\mathbf{z}}
\newcommand{\bZ}{\mathbf{Z}}
\newcommand{\N}{\mathcal{N}}
\newcommand{\diag}{\mathrm{diag}}
\newcommand{\bsigma}{\bm{\sigma}}
\newcommand{\Id}{\mathrm{I}}

\newcommand{\ns}{n_{\mathrm{s}}}
\newcommand{\bdZ}{\bd_{\bZ}}
\newcommand{\bdz}{\bd_{\bz}}

\begin{document}

\title{WaRIance: wavefield reconstruction inversion\\with stochastic variable projection}
\author{G. Rizzuti and T. van Leeuwen}
\date{}

\maketitle

\section{Abstract}

We propose a variation on wavefield reconstruction inversion for seismic inversion, which takes advantage of randomized linear algebra as a way to overcome the typical limitations of conventional inversion techniques. Consequently, we can aim both to robustness towards convergence stagnation and large-sized 3D applications. The central idea hinges on approximating the optimal slack variables involved in wavefield reconstruction inversion via a low-rank stochastic approximation of the wave-equation error covariance. As a result, we obtain a family of inversion methods parameterized by a given model covariance (suited for the problem at hand) and the rank of the related stochastic approximation sketch. The challenges and advantages of our proposal are demonstrated with some numerical experiments.

\newpage
\section{Introduction}

In this paper, we discuss a stochastic approach to wavefield reconstruction inversion \citep[WRI,][]{van2013mitigating}, referred to as WaRIance in the remainder of this paper. WRI was originally conceived an an alternative to full-waveform inversion \citep[FWI: for a summary, refer to][]{virieux2009overview} that relaxes the convergence requirements on the initial model. The adoption of WRI at large has been curbed by the computational demands required to solve the so-called ``augmented'' wave equation, for which direct or iterative solution methods are quite challenging in 3D. The literature around WRI is now extensive, the references most relevant for our work being \cite{wang2016full}, \cite{huang2017full}, \cite{huang2018volume}, \cite{symes2020wavefield}. Only few works had attempted to overcome the computational issues of conventional WRI, for example \cite{wang2016full} and \cite{gholamiWRItime}. Recently, the authors presented a version of WRI where the augmented wavefield is approximated in a sensible manner \citep[for the details, see][]{rizzutiWRI}, albeit with some drawbacks in terms of reconstruction quality and computational overhead. This project can be seen as a way to improve on \cite{rizzutiWRI}. The WaRIance technique here proposed retains the classical well-conditioning properties of WRI in terms of local minimum stagnation, but also can reasonably scale to large 3D problems.

\section{Wavefield reconstruction inversion: a probabilistic perspective}

A probabilistic approach to seismic inversion may consist in modeling additive measurement noise as a Gaussian distribution with a given covariance $\Sd$, that is
\begin{equation}\label{eq:data_model}
    \bd=\R\bu+\bnd,\qquad\bnd\sim\calN(\mathbf{0},\Sd).
\end{equation}
Here, $\bd$ denotes the collected data and $\R$ is the receiver-restriction operator. The unknown $\bu$ represents an underlying wavefield. Contrary to the standard FWI approach, the wave equation can also be formulated in a probabilistic setting:
\begin{equation}\label{eq:weq_model}
    \bq=\Am\bu+\bnq\qquad\bnq\sim\calN(\mathbf{0},\Sq).
\end{equation}
The wave equation operator is denoted by $\Am$ and is parameterized by the physical unknown $\bmm$ (e.g., velocity). The source term is $\bq$ and the source-related covariance is indicated by $\Sq$. Solving for $\bu$ in equation \eqref{eq:weq_model} results in $\bd=\Fm\bq+\bnd+\Fm\bnq$, where $\Fm=\R\Am\inv$ is the ``forward'' operator,  and leads to the following maximum a posteriori problem:
\begin{equation}\label{eq:WRI}
    \min_{\bmm}f(\bmm)=\dfrac{1}{2}\norm{\brm}^2_{\tSd(\bmm)},\qquad\tSd(\bmm)=\Sd+\Fm\Sq\Fm\adj.
\end{equation}
The data residual has been indicated by $\brm=\bd-\Fm\bq$. The weighted norm is defined by
\begin{equation}\label{eq:wnorm}
    \norm{\mathbf{x}}^2_{\Sigma}\defeq\<\Sigma^{-1}\mathbf{x},\mathbf{x}\>,
\end{equation}
for a generic covariance $\Sigma$, where $\<\cdot,\cdot\>$ is the standard Euclidean dot product, and the symbol $^*$ represents the adjoint operation. Note that the formulation in equation \eqref{eq:WRI} is an alternative --- but equivalent! --- derivation of WRI \citep[see also][]{vanleeuwen2019note}. An evident computational hurdle, here, is the need for inverting the perturbed data covariance $\tSd$.

\section{WaRIance: theory}

The goal of this section is to present a practical scheme for the optimization of equation \eqref{eq:WRI}. The starting point is the ``dualization'' \citep{rockafellar1970convex} of the data-misfit term: $1/2\norm{\br}^2_{\Sigma}=\sup_{\by}\<\by,\br\>-1/2\norm{\by}^2_{\Sigma\inv}$. We end up with a saddle-point problem \citep[analogously to][]{rizzutiWRI}:
\begin{equation}\label{eq:WRIlagr}
    \min_{\bmm}\max_{\by}\calL(\bmm,\by)=\<\by,\brm\>-\dfrac{1}{2}\norm{\by}^2_{\tSd(\bmm)\inv},\qquad\tSd(\bmm)=\Sd+\Fm\Sq\Fm\adj.
\end{equation}
Note that the optimal slack variable $\by$ is the solution of the linear system $\tSd(\bmm)\by=\brm$. The gradient of $\calL$ with respect to $\bmm$ is given by $\nabla_{\bmm}\calL=-\Jac(\bmm,\bq+\Sq\Fm\adj\by)\adj\by$. The Jacobian $\Jac(\bmm,\widetilde{\bq})$ is the derivative of the mapping $\bmm\mapsto\Fm\widetilde{\bq}$ with respect to $\bmm$ and the action of its adjoint on $\by$ is the temporal cross-correlation of the (second-order time-derivative of the) forward wavefield $\Fm\widetilde{\bq}$ and the ``backward'' wavefield $\Fm\adj\by$.

A recent discussion in \cite{symes2020wavefield} highlights that a source covariance $\Sq$ that penalizes the equation error away from the point-source position produces a better-conditioned problem for inversion. In \cite{huang2018volume}, for example, $\Sq=\diag(\bsigma^2)$ is a diagonal matrix defined as the inverse of the squared distance from the source location. We exploit the ``source-focusing'' character of this type of covariances by adopting a stochastic low-rank approximation:
\begin{eqnarray}
    \Sq=\diag(\bsigma^2)=\mean_{\bz}\ \bz\bz\adj,\ \ \ \bz\sim\N(\mathbf{0},\diag\bsigma)\quad\implies\quad\calL(\bmm,\by)=\mean_{\bz}\ \calL(\bmm,\by;\bz),\label{eq:covrand1}\\
    \calL(\bmm,\by;\bz)=\<\by,\brm\>-\dfrac{1}{2}\norm{\by}^2_{\tSd(\bmm;\bz)\inv},\quad\tSd(\bmm;\bz)=\Sd+\bdz\bdz{}\adj, \ \bdz=\Fm\bz.\label{eq:covrand2}
\end{eqnarray}
We can then resort to the Monte-Carlo approximation of \eqref{eq:covrand1}:
\begin{equation}\label{eq:WRImc}
    \calL(\bmm,\by)\approx\dfrac{1}{r}\sum_{i=1}^{r}\calL(\bmm,\by;\bz_i)\eqdef\calL(\bmm,\by;\bZ),\qquad\bZ\defeq[\bz_1,\ldots,\bz_{r}]/\sqrt{r},\ \ \bz_i\sim\N(\mathbf{0},\diag\bsigma).
\end{equation}
Here, $\bZ$ is intended as an abstract matrix of rank $r$. This approximation yields the stochastic Lagrangian:
\begin{equation}\label{eq:WRIstochlr}
    \min_{\bmm}\max_{\by}\calL(\bmm,\by;\bZ)=\<\by,\brm\>-\dfrac{1}{2}\norm{\by}^2_{\tSd(\bmm;\bZ)\inv},\qquad\tSd(\bmm;\bZ)=\Sd+\bdZ\bdZ{}\adj, \ \bdZ=\Fm\bZ.
\end{equation}
Note that the optimal multiplier $\tSd(\bmm;\bZ)\by=\brm$ can be calculated via the Woodbury formula
\begin{equation}\label{eq:Woodbury}
    \tSd(\bmm;\bZ)\inv=\Sd\inv-\Sd\inv\bdZ\left[\,\Id+\bdZ{}\adj\,\Sd\inv\bdZ\right]\inv\bdZ{}\adj\,\Sd\inv,
\end{equation}
once the stochastic forward problem $\bdZ=\Fm\bZ$ is solved, assuming that $\Sd$ can be inverted easily. Depending on whether $\bZ$ is drawn at each iteration or kept fixed, one can employ conventional or stochastic optimization for the objective in equation \eqref{eq:WRIstochlr}. The steps required to produce a model update according to the proposed scheme are summarized in Algorithm \ref{alg:WRIstochlr}.
\begin{algorithm}
    \begin{algorithmic}[1]
    \Require{$\bd$, $\Sd$, $\bq$, $\Sq$, $\bmm$}\Comment{Data, data covariance, source, source covariance, physical parameters}
    \Statex
    \State $\bu\leftarrow\Am\inv\bq$,\quad $\br\leftarrow\bd-\R\bu$\Comment{Conventional forward problem: $\mathbf{\ns}$ \textbf{PDE solves}}
    \State $\bz_i\sim\calN(\mathbf{0},\Sq)$ i.i.d.,\quad$\bZ\leftarrow[\bz_1,\ldots,\bz_r]/\sqrt{r}$\Comment{Sampling source covariance approximation}
    \State $\bu_{\bZ}\leftarrow\Am\inv\bZ$,\quad$\bd_{\bZ}\leftarrow\R\bu_{\bZ}$\Comment{Stochastic forward problem: $\mathbf{r}$ \textbf{PDE solves}}
    \Let{$\by$}{$\tSd(\bmm;\bZ)\inv\br$}\Comment{Variable projection, see eq. \eqref{eq:Woodbury}}
    \Let{$\bar{\bu}$}{$\bu+\bu_{\bZ}\,\bd_{\bZ}{}\adj\by$}\Comment{``Augmented'' forward problem}
    \Let{$\mathbf{v}$}{$\Fm\adj\by$}\Comment{Backward problem: $\mathbf{\ns}$ \textbf{PDE solves}}
    \Let{$\Delta\bmm$}{$-\partial_{tt}\bar{\bu}\star\mathbf{v}$}\Comment{Temporal cross-correlation}
    \State \Return{$\Delta\bmm$}
    \end{algorithmic}
    \caption{WaRIance: gradient calculation (for the no-replacement version, fix $\bz_i$'s in line 2)}\label{alg:WRIstochlr}
\end{algorithm}

\section{Numerical experiments}

In this section, we consider a numerical example that is illustrative of the local minimum phenomenon that often spoils FWI. The model consists of a low-velocity Gaussian-shaped lens with a transmission source-receiver configuration with 50 sources and 201 receivers. The starting guess for inversion is an homogeneous model. Data is simulated at 6 Hz and inverted for in the frequency domain. We will compare different objective functionals based on FWI and WRI, and rely on Anderson acceleration as the optimization algorithm \citep{yang2021anderson}.

\begin{figure}[!htb]
  \centering
  \begin{subfigure}[t]{0.23\textwidth}
    \centering
    \includegraphics[width=\textwidth]{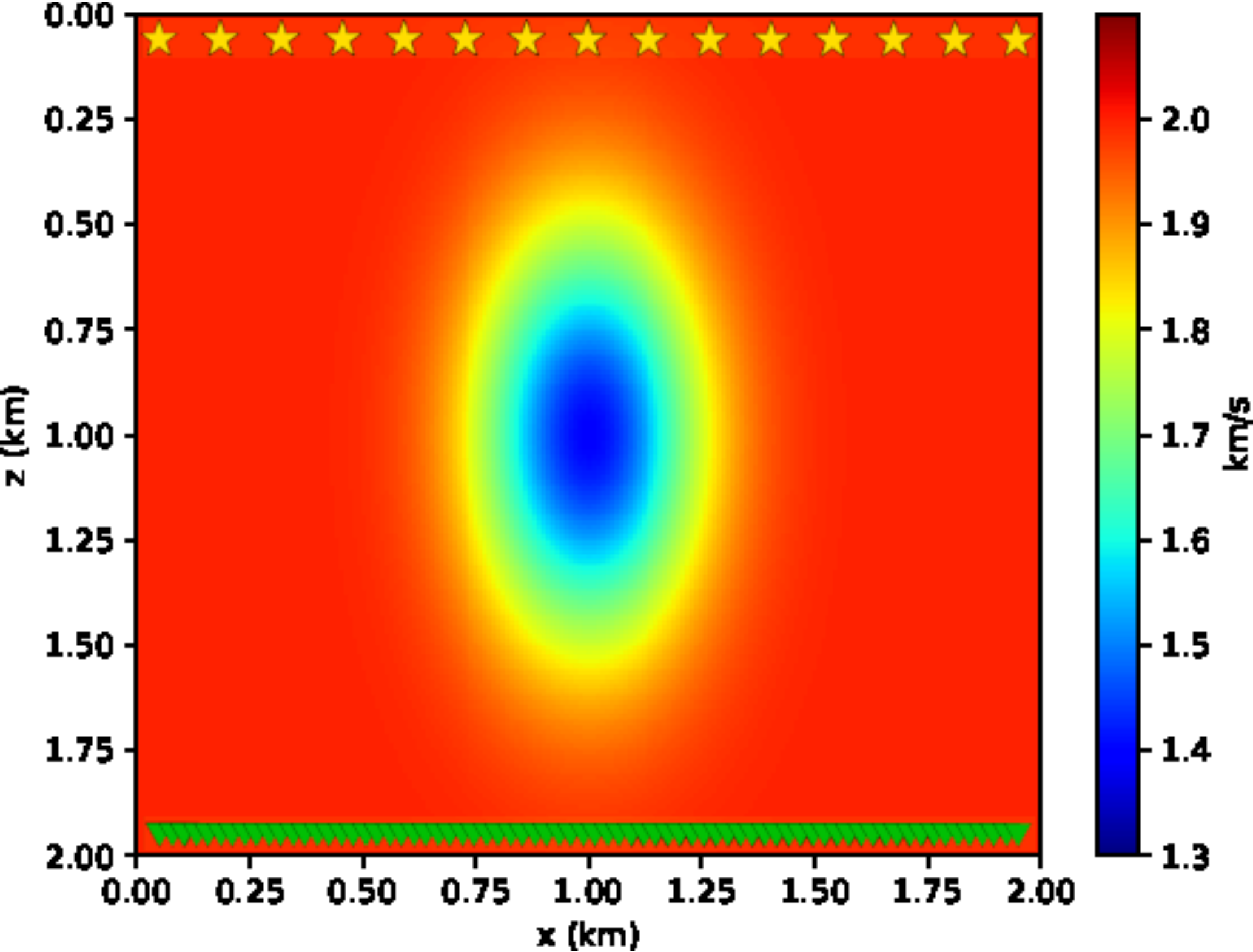}\vspace{-0.25em}
    \caption{True velocity model}\label{fig:GaussianLens_true}
  \end{subfigure}%
  \begin{subfigure}[t]{0.23\textwidth}
    \centering
    \includegraphics[width=\textwidth]{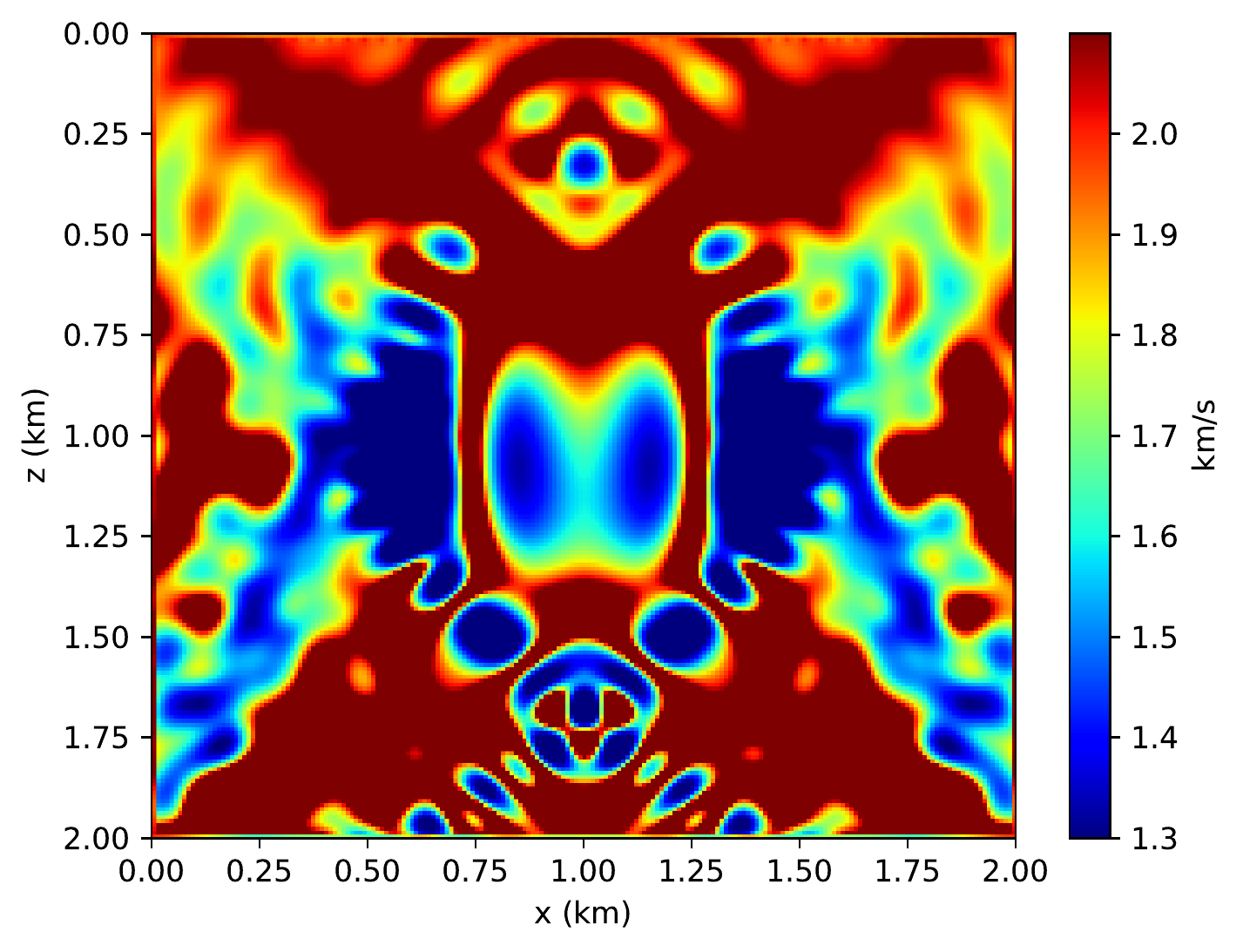}\vspace{-0.25em}
    \caption{FWI}\label{fig:GaussianLens_FWI}
  \end{subfigure}%
  \begin{subfigure}[t]{0.23\textwidth}
    \centering
    \includegraphics[width=\textwidth]{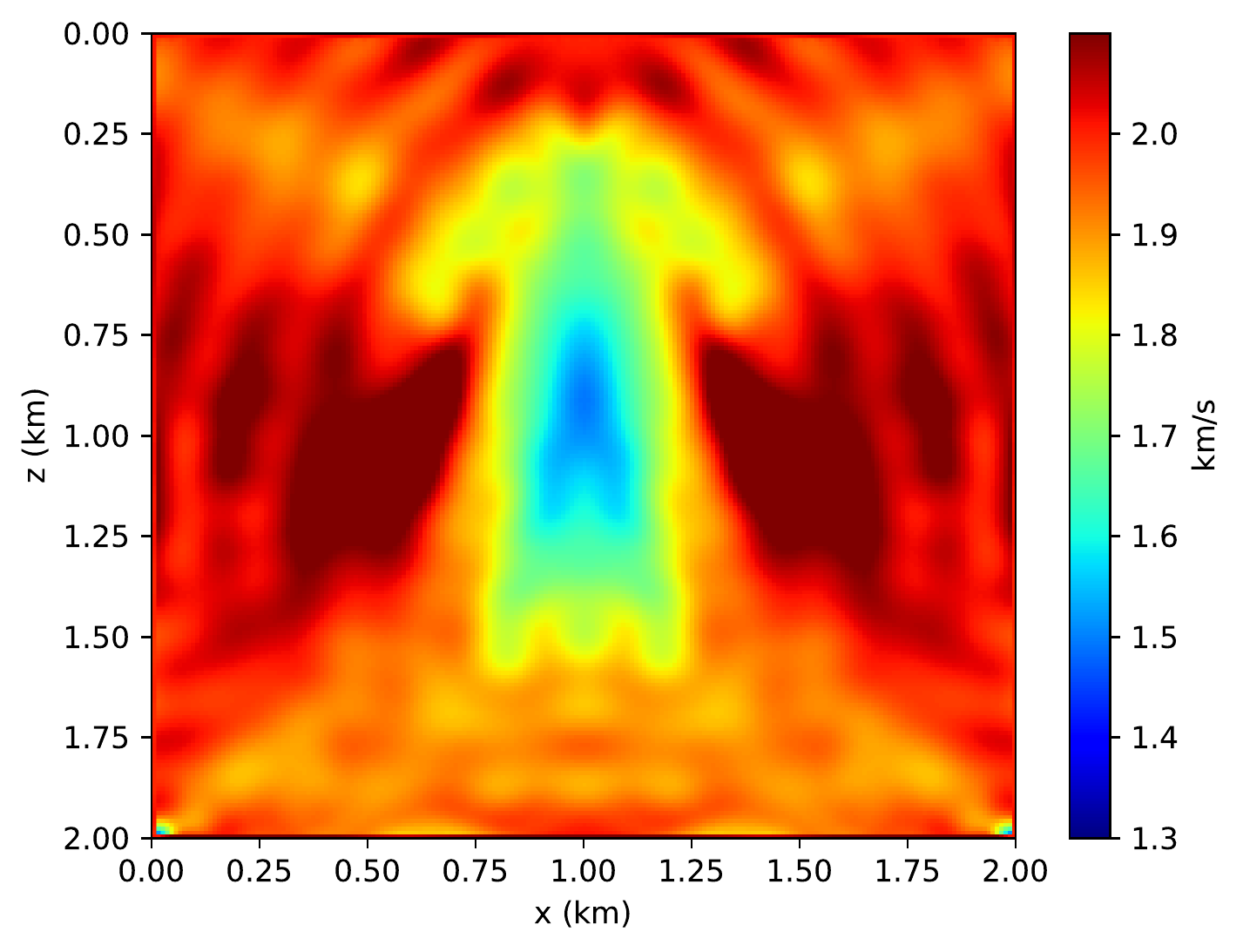}\vspace{-0.25em}
    \caption{WRI}\label{fig:GaussianLens_WRIsrcfcs}
  \end{subfigure}\vspace{-0.5em}%
  \caption{Comparison of inversion results with FWI and source-focusing WRI (equation \ref{eq:WRI}) with source covariance $\Sq$ in equation \eqref{eq:srcfcs}.}\label{fig:GaussianLens}
\end{figure}

\begin{figure}[!htb]
    \centering
    \begin{subfigure}[t]{0.23\textwidth}
        \centering
        \includegraphics[width=\textwidth]{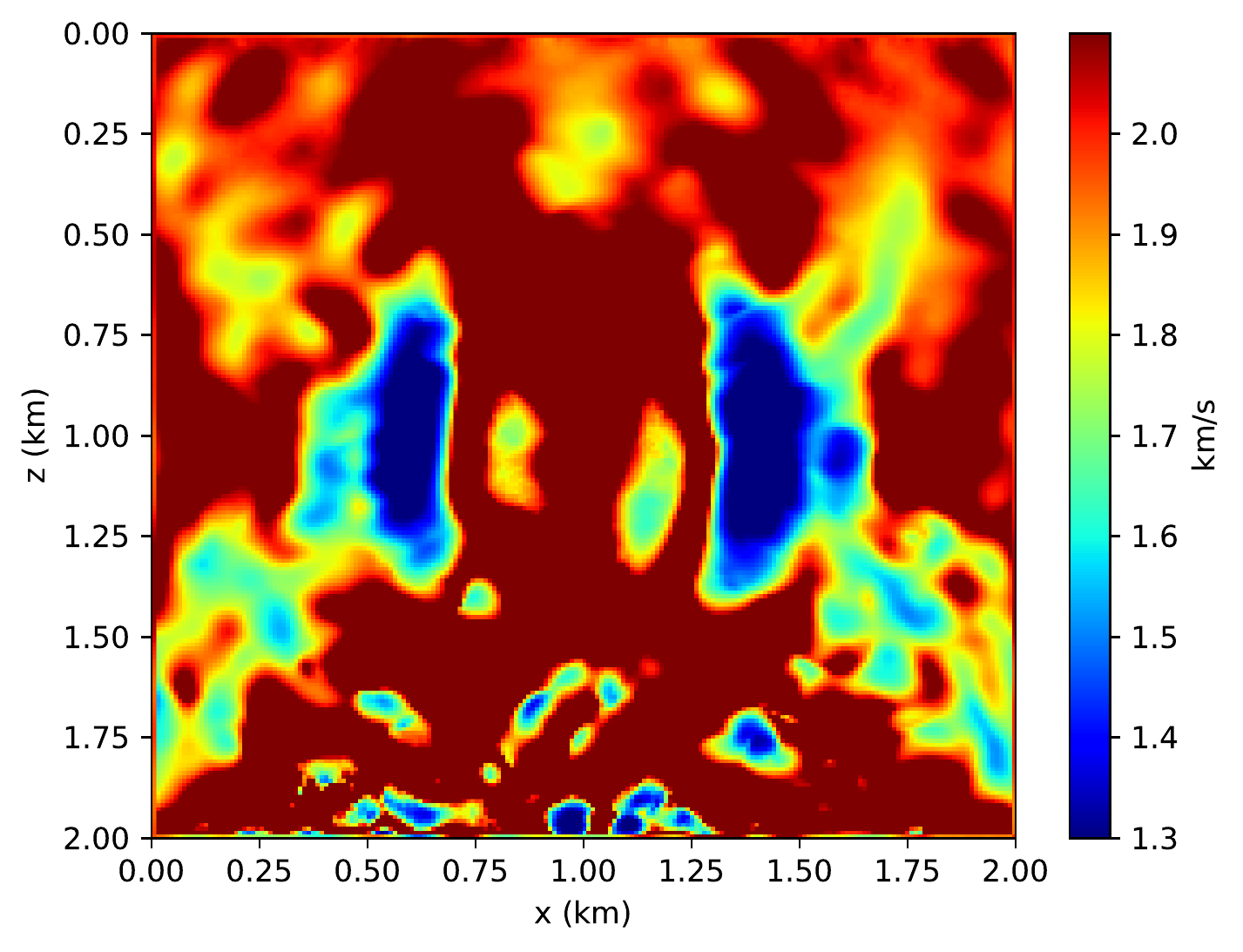}\vspace{-0.25em}
        \caption{$r=\ns$}\label{fig:GaussianLensWRIrandlr_e2_r1}
    \end{subfigure}%
    \begin{subfigure}[t]{0.23\textwidth}
        \centering
        \includegraphics[width=\textwidth]{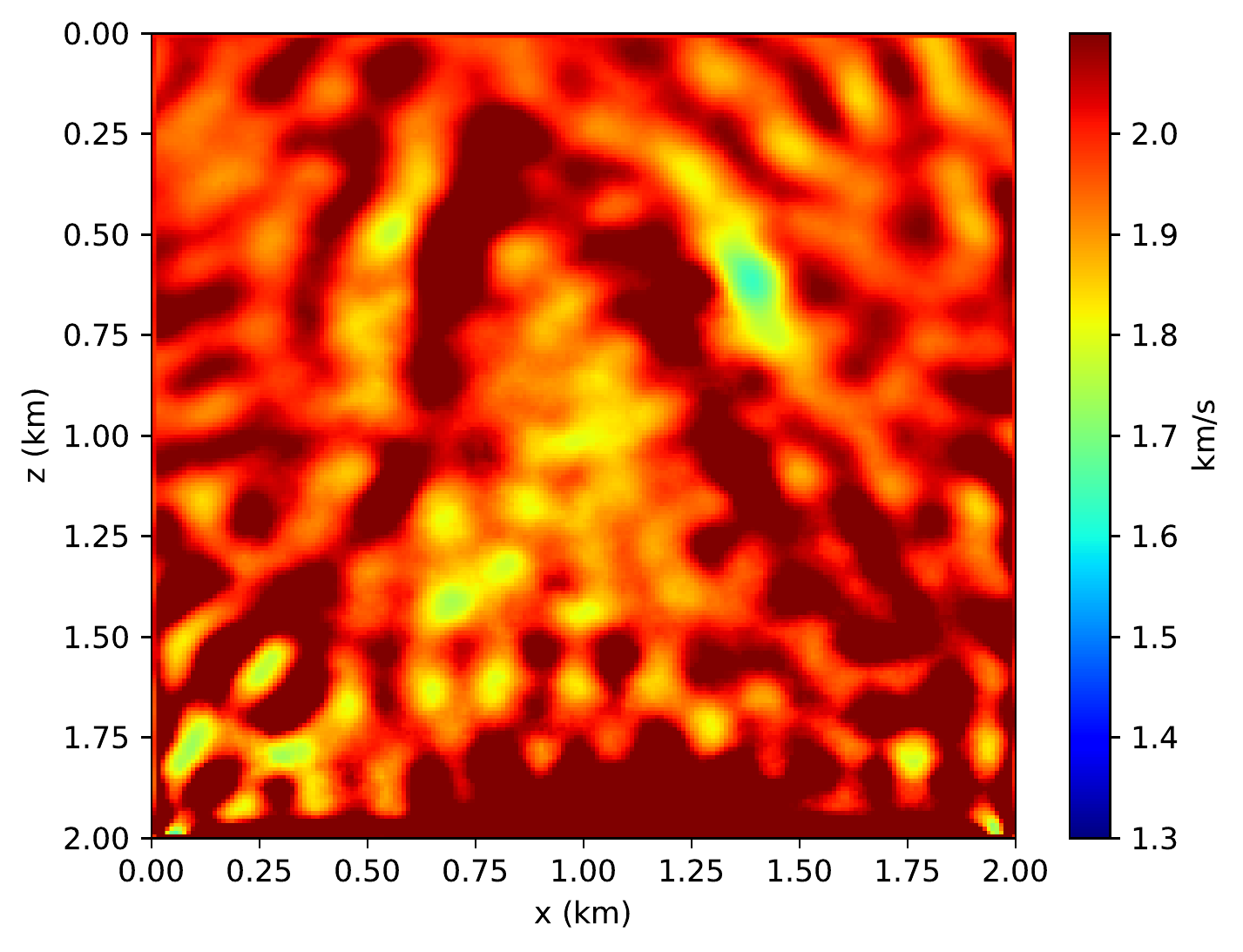}\vspace{-0.25em}
        \caption{$r=10\ns$}\label{fig:GaussianLensWRIrandlr_e2_r10}
    \end{subfigure}%
    \begin{subfigure}[t]{0.23\textwidth}
        \centering
        \includegraphics[width=\textwidth]{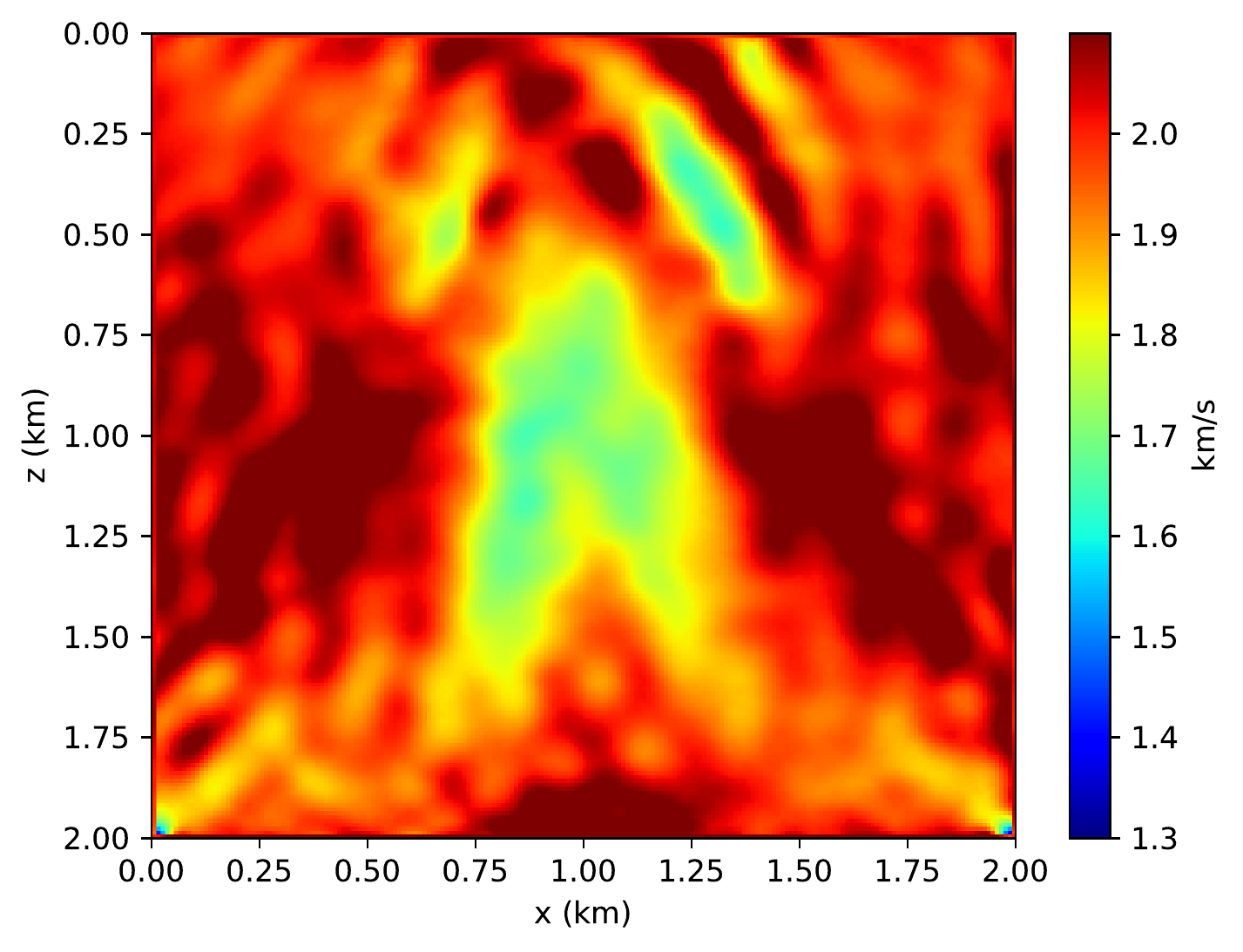}\vspace{-0.25em}
        \caption{$r=30\ns$}\label{fig:GaussianLensWRIrandlr_e2_r30}
    \end{subfigure}%
    \begin{subfigure}[t]{0.23\textwidth}
        \centering
        \includegraphics[width=\textwidth]{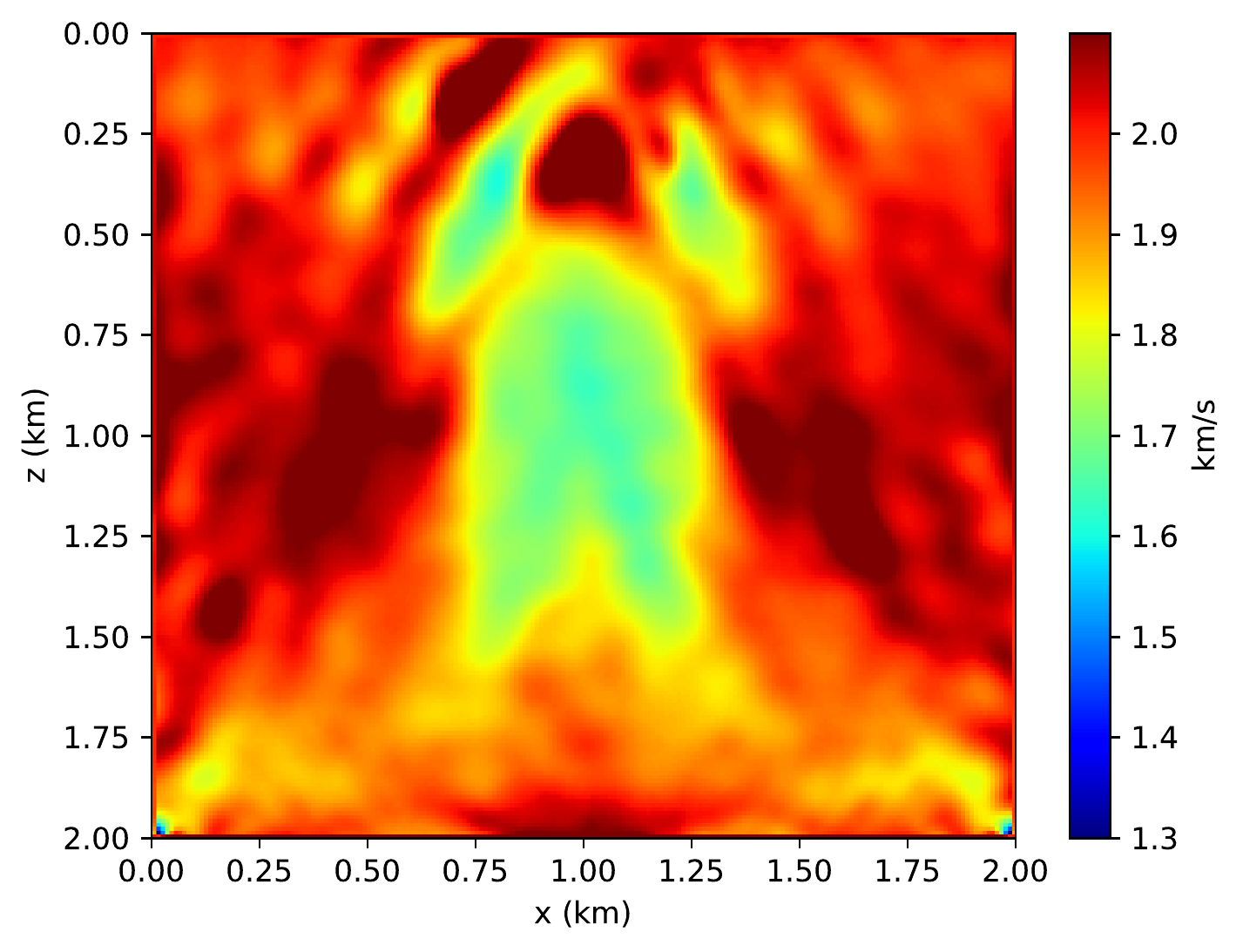}\vspace{-0.25em}
        \caption{$r=50\ns$}\label{fig:GaussianLensWRIrandlr_e2_r50}
    \end{subfigure}\vspace{-0.5em}%
    \caption{Inversion results for WaRIance (equation \ref{eq:WRIstochlr}) with source covariance $\Sq^{\alpha}$ with $\alpha=1$ and different choices of the approximation rank $r$.}\label{fig:GaussianLensWRIrandlr_e2}
\end{figure}

\begin{figure}[!htb]
    \centering
    \begin{subfigure}[t]{0.23\textwidth}
        \centering
        \includegraphics[width=\textwidth]{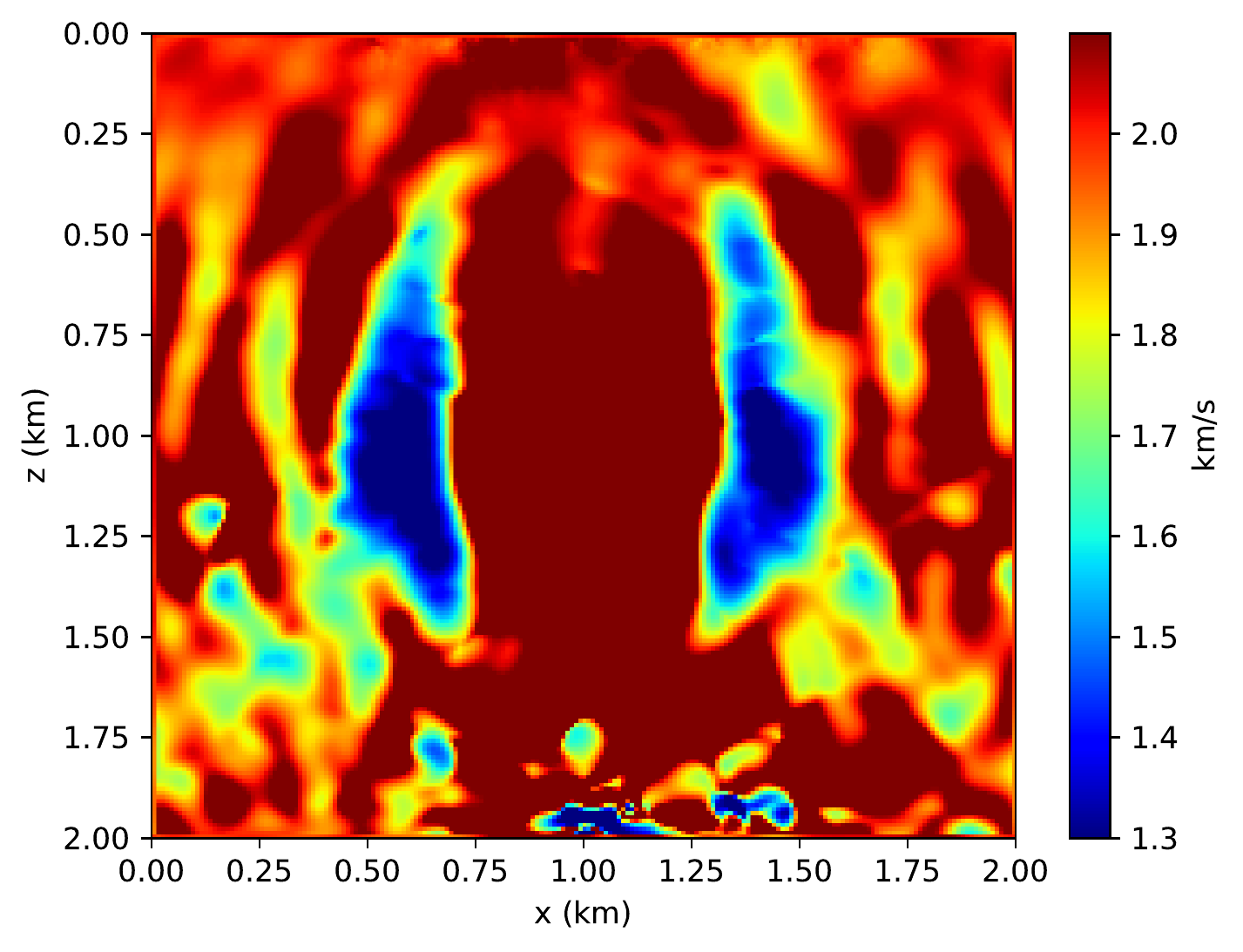}\vspace{-0.25em}
        \caption{$r=\ns$}\label{fig:GaussianLensWRIrandlr_e4_r1}
    \end{subfigure}%
    \begin{subfigure}[t]{0.23\textwidth}
        \centering
        \includegraphics[width=\textwidth]{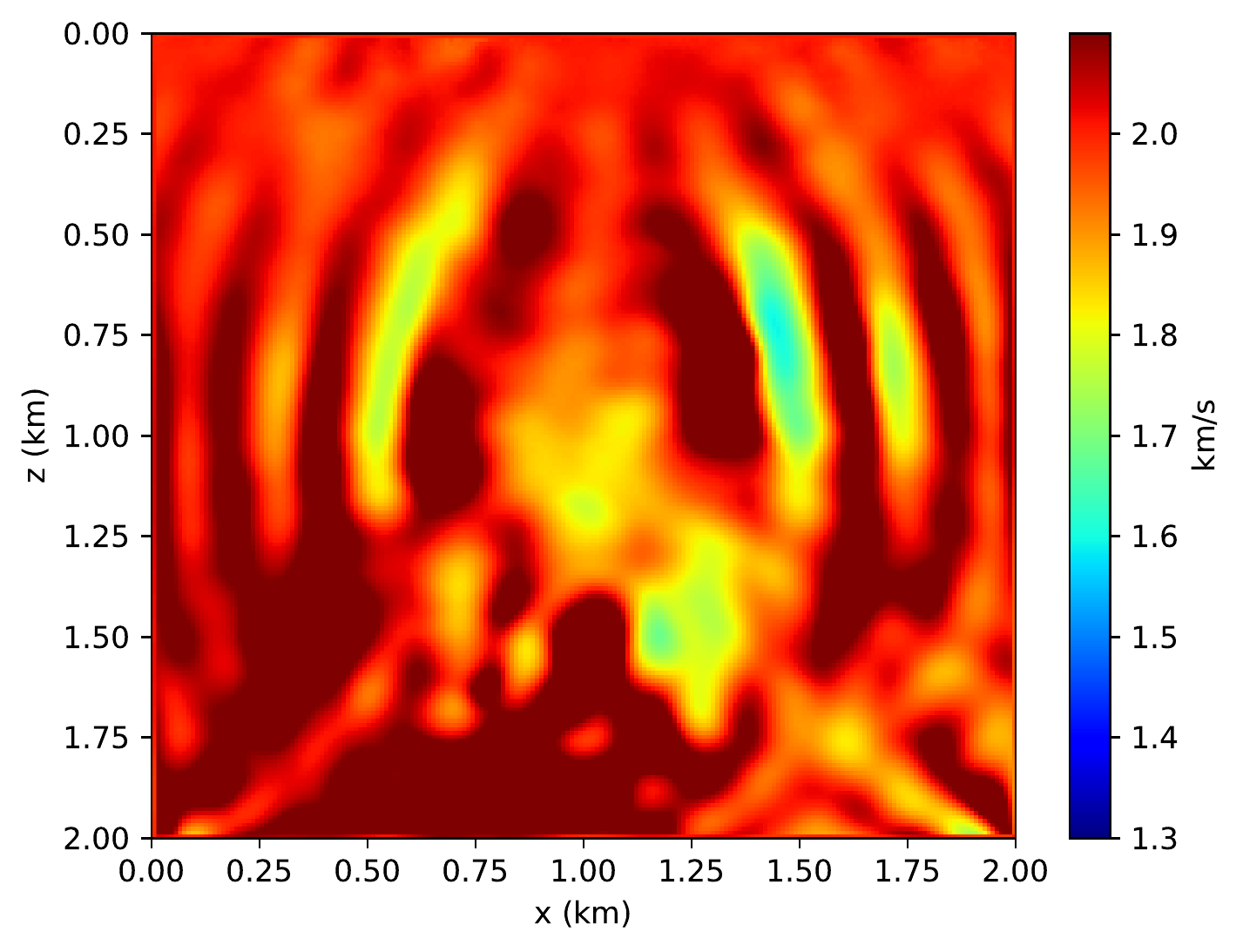}\vspace{-0.25em}
        \caption{$r=2\ns$}\label{fig:GaussianLensWRIrandlr_e4_r10}
    \end{subfigure}%
    \begin{subfigure}[t]{0.23\textwidth}
        \centering
        \includegraphics[width=\textwidth]{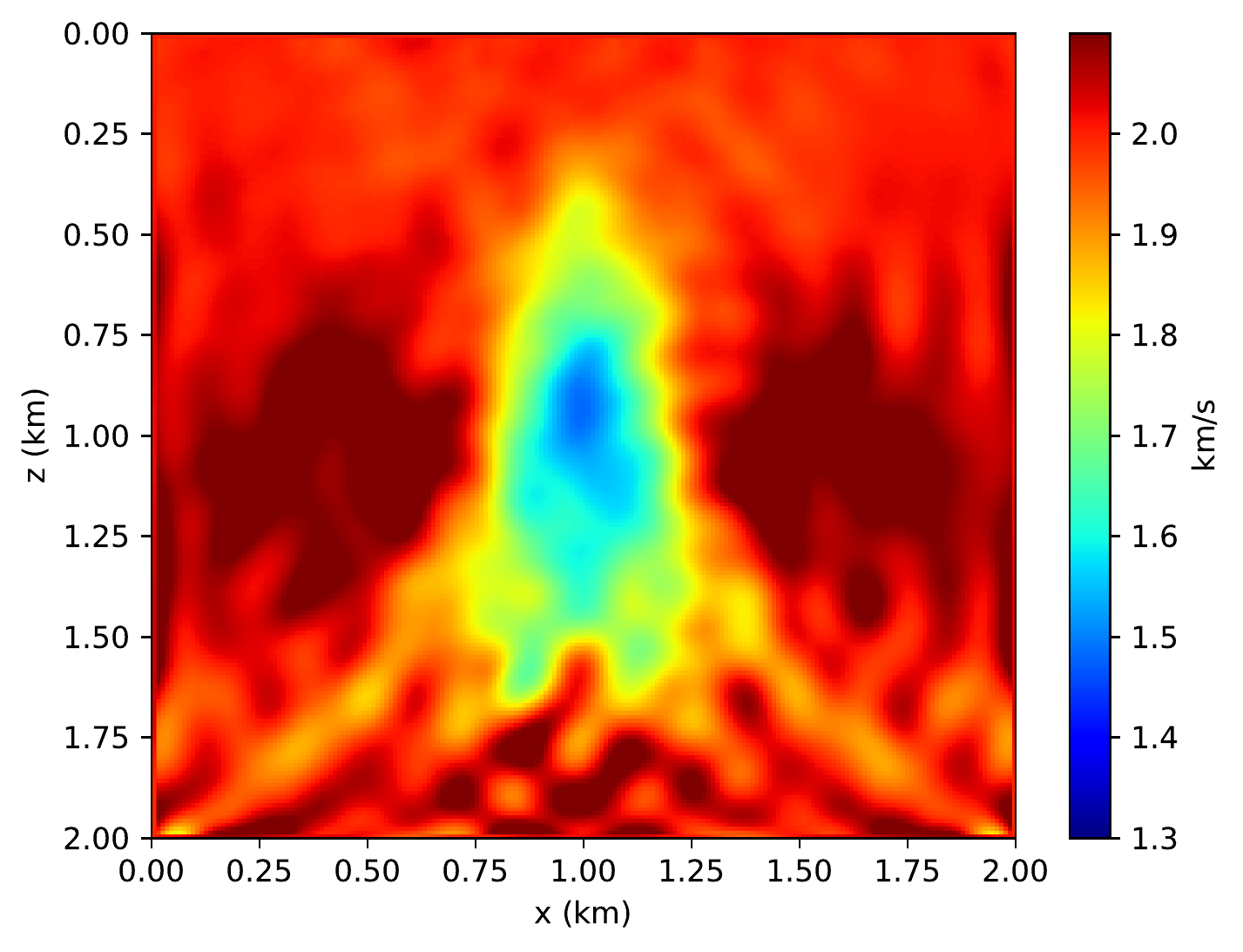}\vspace{-0.25em}
        \caption{$r=3\ns$}\label{fig:GaussianLensWRIrandlr_e4_r30}
    \end{subfigure}%
    \begin{subfigure}[t]{0.23\textwidth}
        \centering
        \includegraphics[width=\textwidth]{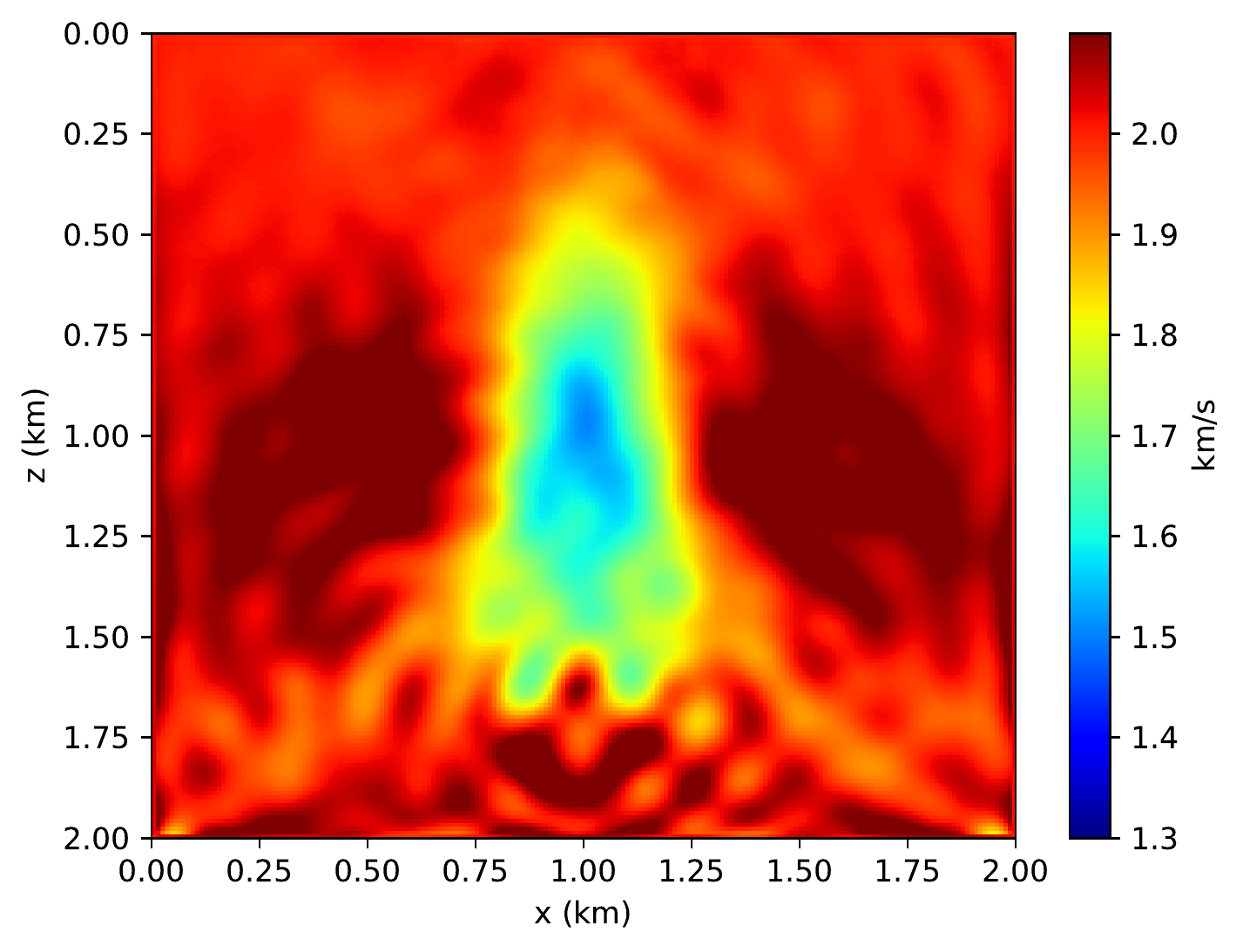}\vspace{-0.25em}
        \caption{$r=4\ns$}\label{fig:GaussianLensWRIrandlr_e4_r50}
    \end{subfigure}\vspace{-0.5em}%
    \caption{Inversion results for WaRIance (equation \ref{eq:WRIstochlr}) with source covariance $\Sq^{\alpha}$ with $\alpha=2$ and different choices of the approximation rank $r$.}\label{fig:GaussianLensWRIrandlr_e4}
\end{figure}

The results for FWI and the source-focusing WRI method of \cite{huang2018volume} are presented in Figure \ref{fig:GaussianLens}, which clearly shows the inability of FWI to converge to the true model. For the source-focusing WRI (equation \ref{eq:WRI}), the source covariance is defined as a diagonal matrix $\Sq=\diag(\bsigma^2)$. The variance function $\bsigma$ is spatially and source-position dependent, e.g.:
\begin{equation}\label{eq:srcfcs}
    \Sq=\diag(\bsigma^2),\quad\bsigma^2(\bx;\bxs)=(\norm{\bx-\bxs}^2+\delta^2)\inv.
\end{equation}
The parameter $\delta$ represents a small regularization term. We test the proposed scheme, WaRIance, for different powers of the covariance $\Sq^{\alpha}$, where the parameter $\alpha$ controls the strength of the source-focusing effect, and for increasing rank of the associated stochastic approximation, as defined in equation \eqref{eq:WRImc}. The results are collected in Figure \ref{fig:GaussianLensWRIrandlr_e2} and Figure \ref{fig:GaussianLensWRIrandlr_e4}. These are actually the mean of different runs obtained by drawing random $\bZ$'s in equation \eqref{eq:WRIstochlr} and kept fixed during the inversion.

\section{Discussion}

The numerical results presented in Figures \ref{fig:GaussianLensWRIrandlr_e2} and \ref{fig:GaussianLensWRIrandlr_e4} demonstrate that WaRIance is a potentially feasible reformulation of WRI that retains its main advantages in terms of well-conditioning, provided that the rank $r$ of the stochastic low-rank approximation is large enough. Generally speaking, for increasing rank, the WaRIance results get closer to deterministic WRI (cf. Figure \ref{fig:GaussianLens_WRIsrcfcs} and Figures \ref{fig:GaussianLensWRIrandlr_e2_r1}--\ref{fig:GaussianLensWRIrandlr_e2_r50}), while choosing the rank too small will make the method converge to a local minimum. The computational costs are roughly determined by the number of PDE solves, which for the proposed scheme is equivalent to $\O(2\ns+r)$, in ``big-o'' notation, as opposed to $\O(2\ns)$ for FWI. However, one should notice that for a source-dependent covariance $\Sq$ (as the one in equation \ref{eq:srcfcs}) the rank $r$ is proportional to $\ns$, e.g. $r=k\ns$, which becomes quite challenging for $k\gg1$. In practice, large-scale applications might have to resort to source encoding \citep{krebs2009fast}, for which a source-independent ``depth-focusing'' $\Sq$ is more apt. In this case, the proposed technique can be quite competitive. More generally, source covariances designed with some source-binning procedure will also help in reducing the costs.

As evidenced by the juxtaposition of the results contained in Figures \ref{fig:GaussianLensWRIrandlr_e2} and \ref{fig:GaussianLensWRIrandlr_e4}, an important theme of WaRIance is related to the ``focusing power'' of the chosen source covariance $\Sq$. Roughly speaking, the more compact the spatial support of a standard deviation function $\bsigma(\bx)$ is, the more accurate its stochastic low-rank approximation is. At the same time, however, the well-conditioning of the problem can degrade, as tested experimentally. Hence, it is important to consider the trade-off between focusing power and well-conditioning.

In the future, we would like to explore stochastic optimization for the objective in equation \eqref{eq:WRIstochlr}, where the random variable $\bZ$ is redrawn at each iteration. We will also study the effects of the specific distribution chosen for $\bZ$. Another interesting research venue is randomized matrix sketching of the perturbed data covariance in equation \eqref{eq:WRI} and its link with WaRIance \citep{tropp2017practical}.

\section{Conclusions}

We proposed a formulation of wavefield reconstruction inversion based on stochastic variable projection: WaRIance. The basic idea is to replace a given model covariance with a low-rank stochastic approximation, with substantial computational benefits that may potentially lead to large-scale applications. Some numerical results show that the proposed technique is competitive against standard methods such as full waveform inversion in terms of robustness towards local search stagnation. The effectiveness of our proposal hinges on the premise that the model covariance can be approximated by a low-rank stochastic sketch. The limitations of these assumptions and further validation will be the subject of future work.

%

\bibliography{biblio}

\end{document}